\documentclass[11pt]{amsart}
\usepackage{graphicx}
\usepackage{amsmath,amssymb}
\newtheorem{thm}{Theorem}[section]
\newtheorem{cor}[thm]{Corollary}
\newtheorem{lem}[thm]{Lemma}

\theoremstyle{definition}

\theoremstyle{remark}
\newtheorem{rem}[thm]{Remark}
\theoremstyle{Example}

\numberwithin{equation}{section}

\begin{document}
\title [ On A  Family Of Elliptic Curves With Positive Ranks  ]
{On A  Family Of Elliptic Curves With Positive Ranks Arising From
Pythagorean Triples}
\author[F.A.Izadi\hspace{10pt}K.Nabardi\hspace{10pt} F.Khoshnam]{F.A.Izadi\hspace{10pt}K.Nabardi\hspace{10pt} F.Khoshnam }
\address{ Mathematics Department Azerbaijan university of  Tarbiat Moallem ,
 Tabriz, Iran  f.izadi@utoronto.ca    farzali.izadi@gmail.com}
\address{ Mathematics Department Azerbaijan university of  Tarbiat Moallem ,
 Tabriz, Iran nabardi@azaruniv.edu}
 \address{ Mathematics Department Azerbaijan university of  Tarbiat Moallem ,
 Tabriz, Iran khoshnam@azaruniv.edu}

\begin{abstract}
The aim of this paper is to introduce a new family of elliptic
curves in the form of  $y^2=x(x-a^2)(x-b^2)$ that have positive
ranks. We first generate a list of pythagorean triples $(a,b,c)$
 and then construct this family of elliptic curves. It turn out that
 this new family have positive ranks and search for the  upper bound for their ranks.
\end{abstract}
\maketitle
{\small{\bf Keywords:} elliptic curves; rank; pythagorean triples\\

\section{Introduction}
An elliptic curve $E$ over a field $F$ is a curve
  that is given by an equation of the form
  \begin{equation}
  Y^2+a_1XY+a_3=X^3+a_2X^2+a_4X+a_6,\quad a_i\in F.
  \end{equation}
 We let $E(F)$ denote the set of points $(x,y)\in F^2$ that satisfy
 this equation, along with a point at infinity denoted O  \cite{4}.\\
 In order for the curve $(1.1)$ to be an elliptic it must be
 smooth, in other words, the three equations
\begin{equation}
 Y^2+a_1XY+a_3Y=X^3+a_2X^2+a_4X+a_6, \\\end{equation}
  $$a_1Y=3X^2+2a_2X+a_4\quad and\quad 2Y+a_1X+a_3=0$$
cannot be simultaneously satisfied by any $(x,y)\in E(\overline{F})$.

If $Char(F)\neq 2$,  we can reduce $(1.1)$ to the following form
 \begin{equation}
Y^2=X^3+aX^2+bX+C
\end{equation}
 with the $discriminant$ :
 \begin{equation}
 D=-4a^3c+a^2b^2+18abc-4b^3-27c^2.
 \end{equation}
 If furthermore,  the $Char(F)$ does not divide $6$, then we get the simplest
 form of
\begin{equation}
 Y^2=X^3+aX+b,
 \end{equation}
 with
 \begin{equation}
 D=-16(4a^3+27b^2).
 \end{equation}
 \begin{rem}
 The elliptic curve is smooth if and only if $ D\neq 0$ \cite{9}.
 \end{rem}

\section{Elliptic curves over $Q$}
 $Mordell$ proved  that on a rational elliptic curve,
 the rational points form a finitely
 generated abelian group, which is denoted by $E(Q)$ \cite{4}.
 Here we can apply the structure theorem for the finitely
 generated abelian groups to $E(Q)$ to obtain a decomposition of
 $E(Q)\cong Z^r\times Tors_E(Q)$,  where $r$ is an integer called the
 $rank$ of $E$ and $Tors_E(Q)$ is the finite abelian group consisting
 of all the elements of finite order in $E(Q)$.\\
 In 1976, ${\it Barry\   Mazur,}$\  proved the following fundamental
 result.
 The torsion group of every elliptic curve is one of the following $15$ cases :
\begin{equation}
\begin{array}{ll}
  \frac{Z}{mZ} & m=1,2,3,...,10,12 \\
   & \\
  \frac{Z}{2Z}\oplus\frac{Z}{mZ} & m=2,4,6,8. \\
\end{array}
\end{equation}
 This shows that there is no points of order $ 11$, and any $n\geq13$.\\
There is an important theorem  proved by $\it Nagell$ and $\it
Lutz$, which tells us how to find all the rational points of
finite order.\\
\begin{thm}(Nagell-Lutz) Let $E$ be given by $y^2=x^3+ax^2+bx+c$
with $a,b,c\in Z$. Let $P=(x,y)\in E(Q)$. Suppose $P$ has finite
order, Then $x,y\in Z$ and either $y=0$ or $y^2|D$.
\end{thm}
\begin{proof}
 ( \cite{8} . $pp\  .\  56$    ).
\end{proof}
\begin{thm}
Let $E$ be given by $y^2=x^3+ax^2+bx+c$ and, $P=(x,y)\in E(Q)$. $P$
has an order $2$ if and only if $y=0$.
\end{thm}
\begin{proof}
( \cite{9}. $ pp\ .77$ ) .
\end{proof}
 On the other hand, it is not known which values of $rank\ r$ are
possible. The current record is an example of elliptic curve over
$Q$ with $rank\geq28$ found by $\it Elkies$ in may 2006 \cite{2}.\\

 In this Paper we first introduce a family of elliptic curves over
$Q$ and show that they have positive rank, then search  for the
largest ranks possible.\vspace{10pt}

\section{ Pythagorean triples}
 A primitive pythagorean triple is a triple of numbers
$(a,b,c)$ so that $a$\ ,\ $b$\ and $c$ have no common divisors and
satisfy
\begin{equation}
a^2+b^2=c^2. \end{equation}
 It's not hard to prove that if one of
$a$ or $b$ is odd then the other is even, then $c$ is always odd.\\
In general , we can generate $(a,b,c)$ by the following relations:
\begin{equation}
\begin{array}{ccc}
  a=i^2-j^2\quad & b=2ij\quad & c=i^2+j^2 \\
\end{array}
\end{equation}
where $gcd\ (i,j)=1$ and $i$, $j$ have oppositive parity.\\
\\
\
The following table gives all possible triples with $i,j<10$.
\\

\begin{center}
\begin{table}[!h]
\begin{tabular}{|c|c|c|c|c|c|}\hline
$i$&$j$&$a=i^2-j^2$&$ b=2ij$&$c=i^2+j^2$&$(a,b,c)$\\\hline
$2$&$1$&$3$&$4$&$5$&$(3,4,5)$\\\hline
$3$&$2$&$5$&$12$&$13$&$(5,12,13)$\\\hline
$4$&$1$&$15$&$8$&$17$&$(15,8,17)$\\\hline
$4$&$3$&$7$&$24$&$25$&$(7,24,25)$\\\hline
$5$&$2$&$21$&$20$&$29$&$(21,20,29)$\\\hline
$5$&$4$&$9$&$40$&$41$&$(9,40,41)$\\\hline
$6$&$1$&$35$&$12$&$37$&$(35,12,37)$\\\hline
$6$&$5$&$11$&$60$&$61$&$(11,60,61)$\\\hline
$7$&$2$&$45$&$28$&$53$&$(45,28,53)$\\\hline
$7$&$4$&$33$&$56$&$65$&$(33,56,65)$\\\hline
$7$&$6$&$13$&$84$&$85$&$(13,84,85)$\\\hline
$8$&$1$&$63$&$16$&$65$&$(63,16,65)$\\\hline
$8$&$3$&$55$&$48$&$73$&$(55,48,73)$\\\hline
$8$&$5$&$39$&$80$&$89$&$(39,80,89)$\\\hline
$8$&$7$&$15$&$80$&$113$&$(15,80,113)$\\\hline
$9$&$2$&$77$&$36$&$85$&$(77,36,85)$\\\hline
$9$&$4$&$65$&$72$&$97$&$(65,72,97)$\\\hline
$9$&$8$&$17$&$144$&$145$&$(17,144,145)$\\\hline
\end{tabular}\\
\ \caption{\footnotesize Generating the primitive  pythagorean triples with $i,j<10$ }
\end{table}
\end{center}

\
\section{Structure Of  The Curves}
\
\
 First we generate a list of primitive pythagorean triples $(a,b,c)$ with
$i,j\leq 1000$. This yields a list of $\ 202461$ triples. Each
$(a,b,c)$ gives rise to the elliptic curve in the form
\begin{equation}
y^2=x(x-a^2)(x-b^2).
\end{equation}
Then we compute the $2-selmer$ $\it ranks$ of these curves as upper
bounds on the $\it Mordell-Weil\  ranks$, finally, by using
$\it Mwrank$, we can obtain the ranks of corresponding curves.\\

\ \section{Relation Between  Euler's Concordant forms and elliptic
curves}
\
\
 In 1780, Euler asked for a classification of those pairs of distinct non-zero integers  M and N
 for which there are integers solutions $(x,y,t,z)$ with $xy\neq0$
 to the system of equation
 \begin{equation}
\begin{array}{c}
  x^2+My^2=t^2\\
 x^2+Ny^2=z^2.
 \end{array}
 \end{equation}
One  can consider Euler's problem  as the problem of the study of the
 elliptic curve over $Q$. i.e :
 \begin{equation}
 E_Q(M,N):\quad y^2=x^3+(M+N)x^2+MNx.
 \end{equation}
 A solution to $(5.1)$ is primitive, if $x,y,t,$ and $z$ are
 positive integers and $gcd(x,y)=1.$
If $E_Q(M,N)$ has positive rank, then there are infinity many
primitive integer solutions to $(5.1)$ \cite{6}. If $E_Q(M,N)$ has
rank $0$, then $(5.1)$ has  a solution if and only if the torsion
group is
$$ \frac{Z}{2Z}\oplus\frac{Z}{8Z}\quad or\quad
\frac{Z}{2Z}\oplus\frac{Z}{6Z}.$$
 We can let the $gcd(M,N)$ be a square-free integer, also we
 can show that $E_Q(M,N)\simeq E_Q(-M,N-M)\simeq
 E_Q(-N,M-N).$ Therefore without loss of generality assume that $M$
 and $N$ are both positive integers.\\ So in $(5.2)$, if we let
 $M=-a^2$ and $N=-b^2$, where $a^2+b^2=c^2$, then
 $$E_Q(-a^2,-b^2):\quad y^2=x^3+(-a^2-b^2)x^2+a^2b^2x=x(x-a^2)(x-b^2)$$
 which is in the form of $(4.1)$. Therefore if we can prove $(4.1)$ has
 a positive rank or has either a torsion group of   $ \frac{Z}{2Z}\oplus\frac{Z}{8Z}\quad or\quad
\frac{Z}{2Z}\oplus\frac{Z}{6Z}$ , then it turns out that in this case
$(5.1)$,  has infinity many solutions.\\
But as we shall see, $(4.1)$ has the torsion group of
$\frac{Z}{2Z}\oplus \frac{Z}{2Z}$. To prove that our family of
elliptic curves has $\frac{Z}{2Z}\oplus\frac{Z}{2Z}$ as a torsion
group,
 among other things, we need to use the following theorem too.

\begin{thm}
The torsion subgroup of $E_Q(M,N)$ are uniquely determined by the following four cases:\\
i) The torsion subgroup of $E_Q(M,N)$ contains $\frac{Z}{2Z}\oplus
\frac {Z}{4Z}$ if $M$ and $N$ are both squares, or $-M$ and $N-M$
are both squares, or if $-N$ and $M-N$ are both squares.\\
\\
ii) The torsion subgroup of $E_Q(M,N)$ is
$\frac{Z}{2Z}\oplus\frac{Z}{8Z}$ if there exists a non-zero integer
$d$ such that $M=d^2u^4$ and $N=d^2v^4$, or  $M=-d^2v^4$ and
$N=d^2(u^4-v^4)$, or $M=d^2(u^4-v^4)$ and $N=-d^2v^4$ where
$(u,v,w)$ forms a pythagorean triple. \\
\\
iii) The torsion subgroup of $E_Q(M,N)$ is $\frac{Z}{2Z}\oplus
\frac{Z}{6Z}$ if there exists integers $a$ and $b$ such that
$\frac{a}{b}\notin\{-2,-1,\frac{-1}{2},0,1\}$ and $M=a^4+2a^3b$ and
$N=b^4+2ab^3$.\\
\\
iv) In all other cases, the torsion subgroup of $E_Q(M,N)$ is
$\frac{Z}{2Z}\oplus\frac{Z}{2Z}$
\end{thm}
\begin{proof}
 ( [6] )
\end{proof}
\section{Results About  The New Family Of Curves}
\
\
 \begin{rem}
  For any pythagorean triple $(a,b,c)$,the elliptic curve in the form $y^2=x(x-a^2)(x-b^2)$
  is smooth. In fact $a\neq b$ and both are nonzero.
 \end{rem}
 \begin{rem}
In  $[3]$, Fouvry and Pomykala lead to an interesting result which
is following. Let $E$ be an elliptic curve in the form
\begin{equation}
y^2=x^3+a(t)x+b(t)
\end{equation} where $a(t) ,\  b(t)\in Z[t]$.
Then the average rank of $E$ is bounded by $2\max\{3deg\ a,2 deg\
b\}.$ Therefore if we change $(4.1)$ to the $(5.1)$ and let one of
the $i$ or $j$ be constant, we would have $a$ and $b$ the
polynomials with degree $8,12$ . So we have
 $r\leq2\max\{3deg a,2degb\}=48$.
\end{rem}
  \begin{lem} The elliptic curve in the form $(4.1)$ has  four points
 of order 2.
 \end{lem}
 \begin{proof}
 It is clear that the points $P_1=(0,0),P_2=(a^2,0),P_3=(b^2,0)$ are
 of order 2. Then $2E(Q)\simeq \frac{Z}{2Z}\oplus \frac{Z}{2Z}.$
 \end{proof}
 \begin{thm}
 Let $E$ be an elliptic curve defined over a field $F$, by the equation
 $y^2=(x-\alpha)(x-\beta)(x-\gamma)=x^3+ax^2+bx+c$, where $Char(F)\neq2$ . For
 $(x',y')\in E(F)$, there  exists $(x,y)\in E(F)$ with
 $2(x,y)=(x',y')$, if and only if $x'-\alpha$, $x'-\beta$, and
 $x'-\gamma$ are squares.
 \end{thm}
 \begin{proof}
 (\cite{4}. Th 4.1. pp.37 ).
 \end{proof}
 \begin{thm}
 The elliptic curve in the form $(4.1)$ doesn't have any point of
 order $4$.
 \end{thm}
\begin{proof}
Let $P=(x,y)\in E(Q)$, such that $4P=O$. Then one of following cases
must be true.
$$2P=(0,0) \qquad or\quad 2P=(a^2,0)\quad or  \qquad 2P=(b^2,0).$$
If $2P=(0,0)$, then $-a^2$ and $-b^2$, are squares, which is a
contradiction. If $2P=(a^2,0)$, then $a^2-b^2$ is a square. So we
have, $a^2-b^2=d^2$ for some $d\in Z$ and $a^2+b^2=c^2$. Therefore
$(\frac{a}{b})^2-1=(\frac{d}{b})^2$ and
$(\frac{a}{b})^2+1=(\frac{c}{b})^2$. It turns out that $1$ is a
congruent number again a contradiction. The case $2P=(b^2,0)$ is
similar.
\end{proof}
\begin{cor}
There is a no point of order $8$ on $(4.1)$ .
\end{cor}
\begin{thm}
The elliptic curve in the form $(4.1)$ does not have any point of
order $6$.
\end{thm}
\begin{proof}
We prove this by theorem $(5.1)$. Let $M=-a^2$ and $N=-b^2$ and
without loss of generality assume that $a^2<b^2$. Because
$E_{Q}(M,N)\simeq E_{Q}(-N,\ M-N)$, we continue the proof with
 $E_{Q}(-N,\ M-N)$ which in this case both of the $-N$ and $M-N$ are
positive integers. Let there exist integers $A$ and $B$ such that
$\frac{A}{B}\notin\{-2,\-1,\ \frac{-1}{2},\ 0,\ 1\}$ and
$-N=b^2=A^4+2A^3B$ and $M-N=b^2-a^2=B^4+2AB^3$. Let $b$ is a even
number, so $A$ is as well and since $b^2-a^2$ is odd, then $B$ must
be odd. Since $gcd\ (a,b)=1$ we have $gcd\ (A,B)=1$. $b^2=A^3(A+2B)$
so $a= t^2$ and $A+2B=s^2$ where $t,s\in Z$. Because $A$ is even, so
$A+2B=s^2$ is as well, thus $A+2B\equiv 0 (mod\ 4)$, in other hand
$A$ is even and square, thus $A\equiv\ 0\ (mod\ 4)$, which  means
that $2\mid B$ which is a contradiction.\\
 Now let $b$ is odd, we conclude that both of $A$ and $B$ are odd.
 So $A+2B=s^2$ is odd and then $s^2\equiv 1 (mod\ 4)$ and
 $a\equiv t^2\equiv 1(mod\ 4)$, then we have $B\equiv 0 (mod\ 2)$,
 which is again a contradiction.
\end{proof}
 \begin{lem}
  For each pythagorean triple $(a,b,c)$,  the elliptic Curve
   $y^2=x(x-a^2)(x-b^2)$ has a positive $\it rank$.
 \end{lem}
 \begin{proof}
 Choose $x=c^2$, then $P=(c^2,\pm abc)$. We show that for each
 $(a,b,c)$, $abc$ does not divide the $discriminant\ D$,  where
  $D=a^4b^4(c^4-4a^2b^2)$.
  If $abc \mid a^4b^4(c^4-4a^2b^2)$ then $c \mid a^3b^3(c^4-4a^2b^2)$.
  Let  $p$ is a prime number such that $p\mid c$ ,  then $p \mid-4a^2b^2$,
  but $c$ is odd, then $p\neq2$ so $p\mid a^2b^2$ and hence $p|a$ or
  $p|b,$ which is a contradiction.
 So $p=(c^2,\pm abc)$  has integer coordinate in which $y=\pm abc$ does not
 divide $D$.
 Therefore by $\it Nagell-Lutz$ theorem $P$ does not have finite order.
 This implies that $r\geq1$.
 \end{proof}
  \begin{cor}
 In the case $M=-a^2$ and $N=-b^2$, where $(a,b,c)$ is a Pythagorean
 triple, the Euler's concordant forms has a infinitely many
 primitive  solution.
 \end{cor}
\section{Numerical Results}
\
\
After searching through $202461$ curves, we found 12 curves with
 ${\it selmer}\  6.$ But none of them  had ${\it
rank}\ 6.$ Also we found $834$ curves with ${\it selmer}\ 5$, leading to $53$ curves of rank $5$.\\
The first curve that generated by first pythagorean triple $(3,4,5)$ has ${\it rank}\ 1.$\\
\\
\ In the following table, we have summarized the results of our
computation.
\begin{center}
\begin{table}[!h]
\begin{tabular}{|c|c|c|}\hline
Rank&number&percent\\\hline
 $rank=1$&$45847$&$22.6$\\\hline
 $rank=2$&$16690$&$8.2$\\\hline
 $rank=3$&$6699$&$3.3$\\\hline
 $rank=4$&$948$&$0.4$\\\hline
 $rank=5$&$53$&$0.02$\\\hline
 $1\leq rank\leq2$&$73204$&$36.1$\\\hline
 $1\leq rank \leq3$&$41381$&$20.4$\\\hline
 $1\leq rank \leq4$&$5906$&$2.9$\\\hline
 $1\leq rank \leq 5$&$384$&$0.1$\\\hline
 $1\leq rank \leq 6$&$2$&$0.0009$\\\hline
 $2\leq rank \leq 3$&$6250$&$3$\\\hline
 $2\leq rank \leq 4$&$4507$&$2.2$\\\hline
 $2\leq rank \leq 5$&$100$&$0.04$\\\hline
 $2\leq rank \leq 6$&$5$&$0.002$\\\hline
 $3\leq rank \leq 4$&$183$&$0.09$\\\hline
 $3\leq rank \leq 5$&$296$&$0.14$\\\hline
 $3\leq rank \leq 6$&$0$&$0$\\\hline
 $4\leq rank \leq 5$&$1$&$0.0004$\\\hline
 $4\leq rank \leq 6$&$5$&$0.002$\\\hline
 $5\leq rank \leq 6$&$0$&$0$\\\hline
\end{tabular}\\
\ \caption{The results of computation.}
\end{table}
\end{center}

 \ In the table $3$, we  have listed the
curves that have selmer equals to $6$, without being able to compute
their exact ranks with MWrank.
 {\footnotesize
\begin{center}
\begin{table}
\begin{tabular}{|c|c|c|c|c|}\hline
i&j&$(a,b,c)$&curve&bound\\\hline
&&&&\\$598$&$53$&$(354795,63388,360413)$&$y^2=x^3-129897530569x^2$&$4\leq
r\leq6$\\&&&$+505788650855590611600x$& \\&&&&\\\hline
&&&&\\$629$&$202$&$(354837,254116,436445)$&$y^2=x^3-190484238025x^2$&$4\leq
r\leq6$\\&&&$+8130585454709316664464x$&\\&&&&\\\hline
&&&&\\$760$&$113$&$(564831,171760,590369)$&$y^2=x^3-348535556161x^2$&$4\leq
r\leq6$\\&&&$+9411982512955600953600x$& \\&&&&\\\hline
&&&&\\$777$&$232$&$(549905,360528,657553)$&$y^2=x^3-432375947809x^2$&$4\leq
r\leq6$\\&&&$+39305500949380532025600x$ &\\&&&&\\\hline
&&&&\\$801$&$560$&$(328001,897120,955201)$&$y^2=x^3-912408950401x^2$&$1\leq
r\leq6$\\&&&$+86586744854271550694400x$&\\&&&&\\\hline
&&&&\\$821$&$242$&$(615477,397364,732605)$&$y^2=x^3-536710086025x^2$&$2\leq
r\leq6$\\&&&$+59813703564011517306384x$&\\&&&&\\\hline
&&&&\\$861$&$788$&$(120377,1356936,1362265)$&$y^2=x^3-1855765930225x^2$&$2\leq
r\leq6$\\&&&$+26681224725077190456384x$&\\&&&&\\\hline
&&&&\\$890$&$457$&$(583251,813460,1000949)$&$y^2=x^3-1001898900601x^2$&$2\leq
r\leq6$\\&&&$+225104091544539413571600x$&\\&&&&\\\hline
&&&&\\$917$&$846$&$(125173,1551564,1556605)$&$y^2=x^3-2423019126025x^2$&$4\leq
r\leq6$\\&&&$+37719046943947124807184x$&\\&&&&\\\hline
&&&&\\$957$&$788$&$(294905,1508232,1536793)$&$y^2=x^3-2361732724849x^2$&$2\leq
 r\leq6$\\&&&$+197833836741502151361600x$&\\&&&&\\\hline
&&&&\\$958$&$691$&$(440283,1323956,1395245)$&$y^2=x^3-1946708610025x^2$&$1\leq
r\leq6$\\&&&$+339790269763746950924304x$&\\&&&&\\\hline
&&&&\\$964$&$173$&$(899367,333544,959225)$&$y^2=x^3-920112600625x^2$&$2\leq
r\leq6$\\&&&$+89987080452485248355904x$&\\&&&&\\\hline
\end{tabular}\\
\ \caption{The curves with selmer-rank $6$.}
\end{table}
\newpage
Table$4$, shows some  curves which  rank 5.
\begin{table}[!h]
\begin{tabular}{|c|c|c|c|c|c|}\hline
n&i&j&$(a,b,c)$& curve&${\it rank}$\\\hline
&&&&&\\1&$65$&$58$&$(861,7540,7589)$&$y^2=x^3-57592921x^2$&$5$\\&&&&$+42145284963600x$&\\&&&&&\\\hline
&&&&&\\2&$206$&$73$&$(37107,30076,47765)$&$y^2=x^3-2281495225x^2$&$5$\\&&&&$+1245523255531937424x$&\\&&&&&\\\hline
&&&&&\\3&$219$&$122$&$(33077,53436,62845)$&$y^2=x^3-3949494025x^2$&$5$\\&&&&$+3124065342026615184x$&\\&&&&&\\\hline
&&&&&\\4&$221$&$74$&$(43365,32708,54317)$&$y^2=x^3-2950336489x^2$&$5$\\&&&&$+2011808689365056400x$&\\&&&&&\\\hline
&&&&&\\5&$226$&$197$&$(12267,89044,89885)$&$y^2=x^3-8079313225x^2$&$5$\\&&&&$+1193125293288351504x$&\\&&&&&\\\hline
&&&&&\\6&$277$&$148$&$(54825,81992,98633)$&$y^2=x^3-9728468689x^2$&$5$\\&&&&$+20206925530689960000x$&\\&&&&&\\\hline
&&&&&\\7&$291$&$130$&$(67781,75660,101581)$&$y^2=x^3-10318699561x^2$&$5$\\&&&&$+26299568174145411600x$&\\&&&&&\\\hline
&&&&&\\8&$298$&$241$&$(30723,143636,146885)$&$y^2=x^3-21575203225x^2$&$5$\\&&&&$+19473940840993453584x$&\\&&&&&\\\hline
&&&&&\\9&$305$&$146$&$(71709,89060,114341))$&$y^2=x^3-13073864281x^2$&$5$\\&&&&$+40786150175724531600x$&\\&&&&&\\\hline
&&&&&\\10&$325$&$132$&$(88201,85800,123049)$&$y^2=x^3-15141056401x^2$&$5$\\&&&&$+57269262954257640000x$&\\&&&&&\\\hline
\end{tabular}\\
\
\caption{Some curves with rank $5$.}
\end{table}
\end{center}
}
 In the following table, we have listed  the independent points of the
curves of table $4$
 {\scriptsize
\begin{center}
\begin{table}
\begin{tabular}{|c|c|}\hline
n&Independent points\\\hline
&\\$1$&$(\frac{57564577194761}{1008016},\frac{29006793653594700125}{1012048064})$,$(\frac{165532287616200}{2745649}
,\frac{505394258095121556600}{4549540393})$\\&\\&$(\frac{6192906993}{64},\frac{311795186829399}{512})$,$(\frac{24834332880}{121},\frac{3321719539155360}{1331})$\\&\\&$(341015696,5742307020800)$\\\hline
&\\$2$&$(\frac{166618634504}{121},\frac{311255416873240}{1331})$,$(\frac{12790926337}{9},\frac{-153963331881884}{27})$\\
&\\&$(1862526649,29434944424380)$,$(\frac{14584697373888197298}{2226990481},\frac{45953060323429949195929519458}{105093907788871})$\\&\\&$
(11173929032 ,1060281679441544)$\\\hline
&\\$3$&$(\frac{1420783000225}{2704},\frac{-3709951931018864055}{140608})$,$
(\frac{3426388189979546}{3150625},\frac{-19862798666292714153406}{5592359375})$\\&\\&$(\frac{3209176809789192}{1100401}
,\frac{20777492819646247103496}{1154320649})$,$(\frac{5079795156916250}{1371241},\frac{145504830321607291308950}{1605723211})$\\&\\&$(11153906082
, 964957876872066)$\\\hline
&\\$4$&$(1883980800,2302931030400)$,$(2049417864,18414019508040)$\\&\\&$(\frac{2442134720068225}{602176}
,\frac{-75833401181142946238625}{467288576})$,$(8778656250,-683241762498750)$\\&\\&$(\frac{389025929026}{9}
,\frac{-234351164774907530}{27})$\\\hline
&\\$5$&$(\frac{40247709912197}{724201},\frac{-3971450274935088970094}{616295051})$,$(\frac{14644921094163784}{1292769}
,\frac{964386979747182474225400}{1469878353})$\\&\\&$(\frac{87950467020096}{6889}
,\frac{504745975500657035040}{571787})$,$(18277955208,1851757920077688)$\\&\\&$(42787752953,7974645953968408)$\\\hline
&\\$6$&$(\frac{52434265914}{249001},\frac{-256293028212914618010}{124251499})$,$(120296250,
-47872494168750)$\\&\\&$(6723284800,3861958531200)$,$(\frac{112595270161250}{16129},
\frac{173400086111756488750}{2048383})$\\&\\&$(\frac{14340640706653}{361},\frac{47589097042950453054}{6859})$\\&\\\hline
&\\$7$&$(\frac{2676650962237850}{1394761},\frac{-230234714875282640110250}{1647212741})$,$(\frac{22163879894522425}{5216656}
,\frac{-554628765666572543285925}{11914842304})$\\&\\&$(\frac{34346962133043282}{5997601}
,\frac{57316484301139284256098}{14688124849})$,$(6253062480,74048765888160)$\\&\\&$(\frac{109261411840568520}{717409}
,\frac{34892314618842917159456520}{607645423})$\\&\\\hline
&\\$8$&$(\frac{730404089870769}{891136},\frac{-37789359740568919672425}{841232384})$,$(\frac{5478549187165109}{6056521}
,\frac{-394874229474026983533710}{14905098181})$\\&\\&$(20665851602,118667705326126)$,$(\frac{73166967363875922}{2745649}
,\frac{9236292756019130201629086}{4549540393})$\\&\\&$(51598853768,8996724544134712)$\\&\\\hline
&\\$9$&$(1837492490,-192369433165070)$,$(2274211682,-192094032181618)$\\&\\&$(\frac{3557867077800}{361}
,\frac{2050506769597435800}{6859})$\\&\\&$(\frac{699532475085000}{32761},\frac{12780541414500071841000}{5929741})$,
$(\frac{831997800678440}{29929},\frac{18315695665342299799960}{5177717})$\\&\\\hline
&\\$10$&$(7819306560,11947900423680)$,$(\frac{947937694496}{121},\frac{18954422023540640}{1331})$\\&\\&$(7908659200
,23645902425600)$,$(\frac{49352010853464722}{4977361},\frac{2582386656676462513905118}{11104492391})$\\&\\&$(\frac{6348468129250}{49}
,\frac{ -15061017382562550750}{343})$\\\hline
\end{tabular}\\
\
 \caption{Independent points of curves of table 3.}
\end{table}
\end{center}
}
\begin{center}
\begin{table}
\begin{tabular}{|c|c|c|c|c|}\hline i&j&$(a,b,c)$& curve&${\it rank}$\\\hline
&&&&\\$26$&$17$&$(387,884,965)$&$y^2=x^3-931225x^2$&$4$\\&&&$+117037883664x$&\\&&&&\\\hline
&&&&\\$43$&$24$&$(1273,2064,2425)$&$y^2=x^3-5880625x^2$&$4$\\&&&$+6903609110784x$&\\&&&&\\\hline
&&&&\\$55$&$34$&$(1869,3740,4181)$&$y^2=x^3-17480761x^2$&$4$\\&&&$+48860938803600x$&\\&&&&\\\hline
&&&&\\$63$&$40$&$(2369,5040,5569)$&$y^2=x^3-31013761x^2$&$4$\\&&&$+142557868857600x$&\\&&&&\\\hline
&&&&\\$66$&$47$&$(2147,6204,6565)$&$y^2=x^3-43099225x^2$&$4$\\&&&$+177422080320144x$&\\&&&&\\\hline
&&&&\\$71$&$58$&$(1677,8236,8405)$&$y^2=x^3-70644025x^2$&$4$\\&&&$+190765045779984x$&\\&&&&\\\hline
&&&&\\$74$&$5$&$(5451,740,5501)$&$y^2=x^3-30261001x^2$&$4$\\&&&$+16271058387600x$&\\&&&&\\\hline
&&&&\\$74$&$23$&$(4947,3404,6005)$&$y^2=x^3-36060025x^2$&$4$\\&&&$+283571724009744$&\\&&&&\\\hline
&&&&\\$74$&$53$&$(2667,7844,8285)$&$y^2=x^3-68641225x^2$&$4$\\&&&$+437644224322704x$&\\&&&&\\\hline
&&&&\\$78$&$35$&$(4859,5460,7309)$&$y^2=x^3-53421481x^2$&$4$\\&&&$+703848328419600x$&\\&&&&\\\hline
\end{tabular}\\
\
 \caption{Some curves with rank $4$.}
\end{table}
\end{center}
\newpage
\begin{center}
\begin{table}[h]
\begin{tabular}{|c|c|c|c|c|}\hline i&j&$(a,b,c)$& curve&${\it
rank}$\\\hline
&&&&\\$13$&$6$&$(133,156,205)$&$y^2=x^3-42025x^2+430479504x$&$3$\\\hline
&&&&\\$13$&$10$&$(69,260,269)$&$y^2=x^3-72361x^2+321843600x$&$3$\\\hline
&&&&\\$19$&$6$&$(325,228,397)$&$y^2=x^3-157609x^2+5490810000x$&$3$\\\hline
&&&&\\$20$&$3$&$(391,120,409)$&$y^2=x^3-167281x^2+2201486400x$&$3$\\\hline
&&&&\\$21$&$8$&$(377,336,505)$&$y^2=x^3-255025x^2+16045795584x$&$3$\\\hline
&&&&\\$21$&$10$&$(341,420,541)$&$y^2=x^3-292681x^2+20511968400x$&$3$\\\hline
&&&&\\$4$&$3$&$(7,24,25)$&$y^2=x^3-625x^2+28224x$&$2$\\\hline
&&&&\\$5$&$2$&$(21,20,29)$&$y^2=x^3-841x^2+176400x$&$2$\\\hline
&&&&\\$7$&$4$&$(33,56,65)$&$y^2=x^3-4225x^2+3415104x$&$2$\\\hline
&&&&\\$8$&$1$&$(63,16,65)$&$y^2=x^3-4225x^2+1016064x$&$2$\\\hline
&&&&\\$9$&$2$&$(77,36,85)$&$y^2=x^3-7225x^2+7683984x$&$2$\\\hline
&&&&\\$2$&$1$&$(3,4,5)$&$y^2-25x^2+144x$&$1$\\\hline
&&&&\\$3$&$2$&$(5,12,13)$&$y^2=x^3-169x^2+3600x$&$1$\\\hline
&&&&\\$4$&$1$&$(15,8,17)$&$y^2=x^2-289x^2+14400x$&$1$\\\hline
&&&&\\$5$&$4$&$(9,40,41)$&$y^2=x^3-1681x^2+129600x$&$1$\\\hline
&&&&\\$6$&$1$&$(35,12,37)$&$y^2=x^3-1369x^2+176400x$&$1$\\\hline
\end{tabular}\\
\ \caption{Some curves with rank 3,2, and 1.}
\end{table}
\end{center}
}

\end{document}